%% file: main.tex
\documentclass[letterpaper,11pt]{article}
\usepackage{amsmath, amssymb, amsthm, mathtools}
\usepackage[margin=1in]{geometry}
\usepackage{graphicx}
\usepackage[colorlinks=true, urlcolor=blue, citecolor=black, linkcolor=black]{hyperref}
\usepackage{cleveref}
\usepackage{subcaption}
\usepackage{adjustbox} % for padding around figures

\newtheorem{theorem}{Theorem}
\newtheorem{definition}[theorem]{Definition}
\newtheorem{example}[theorem]{Example}

\DeclareMathOperator{\Aut}{Aut}

\begin{document}
\title{Building with Blocks: Enumerating Polyforms on Tilings}
\author{
    Bert Dobbelaere\textsuperscript{1},
    Peter Kagey\textsuperscript{2},
    Drake Thomas\textsuperscript{3},
    and
    Andr\'es R. Vindas-Mel\'endez\textsuperscript{4}
    \vspace{10pt}\\
    \textsuperscript{1}Redwire Europe; bert.dobbelaere@redwirespace.eu\\
    \textsuperscript{2}Department of Mathematics and Statistics, Cal Poly Pomona; pkagey@cpp.edu\\
    \textsuperscript{3}Independent; grahamsnumberisbig@gmail.com\\
    \textsuperscript{4}Harvey Mudd College; avindasmelendez@g.hmc.edu}
\date{}
\maketitle
\thispagestyle{empty}

\begin{abstract}
  In areas as diverse as contemporary art, play structures, climbing equipment, and modular construction toys, we see the presence of building block-like polyhedral complexes, which are generalizations of the pieces in the game Tetris.
  We give an algorithm for counting the number of $n$-celled structures on polygonal and polyhedral cells of certain periodic two- and three-dimensional tilings; moreover, we count these structures up to translations, rotations, and reflections of the tiling.
  We describe this algorithm with respect to structures in the snub square tiling, provide numerical data related to existing three-dimensional art and structures, and suggest puzzles based on these constructions.
\end{abstract}

\section*{Introduction}
In Jiangmei Wu's 2017 large-scale installation \emph{Synergia}, shown in \Cref{subfig:JiangmeiWu} and described in the 2018 Bridges Conference Proceedings \cite{WuInchbald2018}, hundreds of copies of the space-filling bisymmetric hendecahedra are arranged to form an illuminated structure in front of Eero Saarinen's North Christian Church.
This is just one of many examples of \emph{polyhedral complexes} arising in contemporary art, architecture, toy design, science communication, and beyond.
We also see such complexes in Tomás Saraceno's 2012 \emph{Cloud City}, installed on the roof garden of the Metropolitan Museum of Art, which suspends a walkable network of polyhedral cells in space;
Richard Dattner's modular playground structures based on cuboctahedra attached face-to-face (\Cref{fig:DatterPlaygroundCuboctahedra}); Robin Houston's LEGO-like building blocks built from tetragonal disphenoids (\Cref{fig:disphenoidTetrahedral}); the ZIGZAG modular climbing wall composed of tetrahedral and octahedral cells (\Cref{fig:octtet}); and Bih-Yaw Jin and Chia-Chin Tsoo's tubular glass model of the tetrahedral-octahedral structure present in perovskite crystal, which appeared in the 2015 Bridges Conference Art Exhibition \cite{Jin2015}.

Simultaneously, there is an artistic and mathematical impulse to \textit{enumerate everything}, seen in works ranging from conceptual artist Sol Lewitt's series of sculptures ``Incomplete Open Cubes'' --- which enumerates all possible ways of \emph{not} making a cube --- to Jorge Luis Borges's speculative fiction short story ``The Library of Babel'' about a library that contains all possible books ``of four hundred and ten pages''.
Here we combine these two ideas to ask the question: given a collection of building blocks subject to some constraints, \emph{how do we enumerate all possible structures}?

\input{figure/bisymmetric_hendecahedron}

Still, there are more examples of polyhedral structures arising in contemporary art and mathematical outreach.
In Olafur Eliasson's 2024 \textit{Rehearsal room for spatial imagination} at The Geffen Contemporary at MOCA, we see tube lights reflecting in a human-scale, minimal octahedral kaleidoscope in a corner of a room to form a polyhedral complex of rhombic dodecahedra with full octahedral symmetry;
Antony Gormley's 55-foot-high, 1700-pound sculpture \emph{Chord}, which consists of 33 polyhedra, hangs between MIT's Mathematics and Chemistry Departments; and
Conrad Shawcross's sculpture ``The Dappled Light of the Sun'' consists of thousands of steel tetrahedra welded together, with a total weight of thirty tons.

Restricting attention to polyhedral complexes consisting only of truncated octahedra, we see in \Cref{fig:TOHYBYCOZO},
HYBYCOZO's $5$-cell truncated octahedral installation at the Desert Botanical Garden in Phoenix, AZ;
in \Cref{fig:TOJasper} Jasper Bown's $10$-celled ceramic sculpture \emph{Polytopes amidst a Permutahedral Tiling}, composed of ten truncated octahedra and exhibited at the Bridges Art Exhibition at the 2025 Joint Mathematics Meetings \cite{Bown2024}; and
in \Cref{fig:TOWARD} Studio Infinity's $38$-celled installation illustrating that ``Truncated Octahedra Work As Rhombic Dodecahedra.''

\input{figure/truncated_octahedra}

\subsection*{Generalized polyforms}
Perhaps the earliest of these structures that humans attempted to enumerate are the $2$-dimensional \textit{polyominoes}, the edge-connected components of the square grid, which were named by Sol Golomb in 1953, popularized for mathematical audiences in Martin Gardner's \emph{Scientific American} ``Mathematical Games'' column in 1960, and entered the mainstream imagination through the popularization of Tetris in the mid-1980s.

Since the 1960s, mathematicians and puzzle designers have given names to components of various tilings of $\mathbb{R}^2$, including but not limited to polyiamonds (triangular), polyhexes (hexagonal), polyocts (truncated square), polykagomes (trihexagonal), polybirds (rhombitrihexagonal), polycairos (Cairo), and polykites (deltoidal trihexagonal tiling).
There are also named $3$-dimensional, face-connected analogs, such as polycubes (cubic), polyrhombs (rhombic dodecahedral), and splatts (truncated octahedral).
Many of these polyforms have been cataloged by George Sicherman \cite{Sicherman}, and some appear in Kadon Enterprises puzzles, such as the ``Ochominoes'' and ``Birds And Bees'' puzzles, the former of which appeared in Kate Jones's printed posterboard piece at the 2021 Joint Mathematics Meetings Bridges Art Exhibition \cite{Jones}.
In general, all such objects are considered \emph{polyforms}.
When identified up to translation, they are called \emph{fixed polyforms}, and when identified further up to rotations and reflections, they are called \emph{free polyforms}.

We give a general framework for enumerating all of the above examples, as well as infinitely many other families of polyforms, which live on \textit{periodic} tessellations of $n$-dimensional Euclidean space, which are called \textit{tilings} in $\mathbb{R}^2$ and \textit{honeycombs} in $\mathbb{R}^3$.
Using this framework, we have confirmed existing counts of polyforms and described over $30$ new families, recording these sequences on the On-Line Encyclopedia of Integer Sequences (OEIS) \cite{OEIS}.
One of the contributions (sequence A385024) counts the polyforms that can be built with Robin Houston's 3D-printed tetragonal disphenoidal building-block toy, shown in \Cref{fig:disphenoidTetrahedral}.

\input{figure/disphenoid_tetrahedron}

\subsection*{Preliminaries}
When counting polyforms on a tessellation of $\mathbb{R}^d$ by $d$-dimensional polytopal (e.g. polygonal or polyhedral) cells, we begin by constructing the dual graph with each vertex of the graph corresponding to a cell in the tessellation, and each edge in the graph corresponding to a cell-to-cell connection along a $(d-1)$-dimensional facet (e.g. edge or face). This allows us to describe polyforms in graph-theoretic terms.
\begin{definition}
  A \textbf{periodic Euclidean graph} is an embedding of a graph into $\mathbb{R}^d$ for which there exist $d$ linearly independent vectors whose translations are graph symmetries (automorphisms).
  The \textbf{translation subgroup} is the subgroup of the automorphism group consisting only of translations and is denoted $T_G \leq \Aut(G)$.
\end{definition}

For each tessellation, there are various families of polyforms, depending on which symmetries (viewed as subgroups of the automorphism group) are considered. Because we are enumerating polyforms by the number of cells, we call an $n$-celled polyform simply an $n$-polyform.

\begin{definition}
  A \textbf{free polyform} on a graph $G$ is an equivalence class of induced subgraphs, where two subgraphs are equivalent if one can be mapped to the other by an automorphism of $G$.
  A \textbf{fixed polyform} is an equivalence class of induced subgraphs taken up to the translation subgroup.

\end{definition}
As a technical note, we only count \textit{quasi-transitive} graphs, which are graphs with a finite number of vertex orbits under the action of the automorphism group. Also, our definition implicitly depends on an embedding of $G$ into $\mathbb{R}^n$, where we only consider automorphisms of $G$ that are isometries of the embedding.

\section*{Overview of technique}
Attempting to enumerate polyforms by hand in an \textit{ad hoc} way---especially in three dimensions---can be a subtly difficult problem. Natasha Rozhkovskaya and Michael Reb demonstrate that in one of Sol Lewitt's enumerations, he omits one of the $122$ of ``Incomplete Open Cubes,'' and duplicates another \cite{Rozhkovskaya2015}.
Similarly, in their 2023 Bridges paper, Robert Bosch and Ilana McNamara demonstrate that in LeWitt's ``All Three-Part Variations on Three Different Kinds of Cubes,'' in which he missed one of the $57$ possible variations \cite{Bosch2023}.
To avoid omitting or duplicating any polyforms, we propose an algorithm to ensure that all possible free polyforms for a specified tessellation are counted precisely once.

The algorithm begins by specifying the orbits of vertices in the underlying dual graph under the automorphism group of the embedding of the graph, and then inductively grows the set of $n$-polyforms from the $(n-1)$-polyforms by adding each possible adjacent cell to each possible $(n-1)$-polyform, and then deduplicating by choosing a canonical representative for each orbit.
The resulting larger polyform is called a \textit{child} of the smaller polyform.
We illustrate this algorithm using a running example that counts the free polyforms on the snub trihexagonal tiling,  shown in \Cref{fig:SnubTrihexagonal}, which is directly analogous to the procedure for counting three-dimensional polyforms.

\input{figure/snub_tri_dual_graph}

\begin{example}
  We begin by embedding the dual graph of the snub trihexagonal tiling in the plane.
  We place a vertex at $(0,0)$ to represent the center of a hexagonal tile, such that translations by $\vec{b}_1 = (5,\sqrt{3})$ and $\vec{b}_2 = (1,3 \sqrt{3})$ correspond to isometries of the tiling.
  These vectors are linearly independent and serve as a basis, allowing us to describe the translation subgroup $T_G < \Aut(G)$ as integer linear combinations of these translations.
  This basis is convenient because all other vertices have rational coordinates with respect to it, and we will describe all vertices and edges in terms of this basis.

  The automorphism group of the snub trihexagonal tiling is the wallpaper group $p6$, which is generated by the following affine transformations acting on homogeneous coordinates.
  The following three transformations correspond, respectively, to translations in the $\vec{b}_1$ and $\vec{b}_2$ directions and $60^\circ$ rotations around the origin:
  \[\Aut(G) = \left\langle
    \left[
      \begin{array}{rr|r}
        1 & 0 & 1\\
        0 & 1 & 0\\\hline
        0 & 0 & 1
      \end{array}
    \right],
    \left[
      \begin{array}{rr|r}
        1 & 0 & 0\\
        0 & 1 & 1\\\hline
        0 & 0 & 1
      \end{array}
    \right],
    \left[
      \begin{array}{rr|r}
        1 & 1 & 0\\
        -1 & 0 & 0\\\hline
        0 & 0 & 1
      \end{array}
    \right]
    \right\rangle\!.
  \]
\end{example}

Next, we describe the vertices and edges of the graph by specifying a representative for each vertex orbit under the automorphism group, choosing the lexicographically earliest vertex with nonnegative coordinates. Note that these vertex orbits correspond precisely to the $1$-polyforms.

\begin{example}
  There are three orbits of vertices under the automorphism group, which correspond to the hexagonal cells, the triangular cells that are not adjacent to a hexagonal cell, and the triangular cells that are adjacent to a hexagonal cell.
  These are shown in \Cref{fig:snubTri1}, and coordinates
  \(v_1 = \left(0,0\right)\),
  \(v_2 = \left(\frac{1}{3},\frac{1}{3}\right)\), and
  \(v_3 = \left(\frac{2}{21},\frac{11}{21}\right)\), respectively.

\input{figure/snub_tri1}
\end{example}

Next, we describe the edges of the graph by specifying the edges that are incident to each of the vertex orbit representatives.

\begin{example}
  For the vertex orbit representative $v_1$, $v_2$, and $v_3$ of the snub trihexagonal tiling, we specify the adjacent vertices (with respect to the translation basis): \begin{align*}
    E_1 &= \left\{
      \left(\frac{-10}{21},\frac{8}{21}\right),
      \left(\frac{-8}{21},\frac{-2}{21}\right),
      \left(\frac{2}{21},\frac{-10}{21}\right),
      \left(\frac{10}{21},\frac{-8}{21}\right),
      \left(\frac{8}{21},\frac{2}{21}\right),
      \left(\frac{-2}{21},\frac{10}{21}\right)
    \right\}
    \\
    E_2 &= \left\{
      \left(\frac{8}{21},\frac{2}{21}\right),
      \left(\frac{11}{21},\frac{8}{21}\right),
      \left(\frac{2}{21},\frac{11}{21}\right)
    \right\}
    \\
    E_3 &= \left\{
      \left(0,0\right),
      \left(\frac{1}{3},\frac{1}{3}\right),
      \left(\frac{13}{21},\frac{-2}{21}\right)
    \right\}
  \end{align*}
\end{example}

After we have described the graph $G$, we need to be able to tell if two induced subgraphs correspond to the same polyform, which we do by specifying a canonical representative for each polyform.

\begin{definition}
    The \textbf{canonical representative} for a polyform up to $R$, where $R$ is a subgroup of $\Aut(G)$ which contains the translation group $T_G$, is the induced subgraph whose vertex set is nonnegative and which is lexicographically minimal when the vertices are sorted lexicographically.
\end{definition}

We compute this canonical representative by applying all possible \emph{orientations}, and then shifting each result so that the coordinates are nonnegative and as small as possible.

\begin{definition}
  Given a graph $G$ with automorphism group $\Aut(G)$ and translation subgroup $T_G$, we say that the quotient group $\Aut(G)/T_G$ partitions the \textbf{orientations} of $G$.
\end{definition}

\begin{example}
  In the case of the snub trihexagonal tiling, the orientations are $\Aut(G)/T_G \cong C_6$, which correspond to the six rotations of the hexagon.
  For coset representatives, we can choose matrices of the form
  \[ A^n = \left[
      \begin{array}{rr|r}
        1 & 1 & 0\\
        -1 & 0 & 0\\\hline
        0 & 0 & 1
      \end{array}
    \right]^n
  \] for $n = 0, 1, \dots, 5$.

  Next we compute the canonical representation for $P = \{(-10/21,8/21), (0,0)\}$ which results from extending the (hexagonal) $1$-polyform $\{(0,0)\}$ representative to one of its neighboring vertices.
  We apply all six orientation transformations and translate to the first quadrant to get the following candidates.
  \begin{alignat*}{2}
    A^0 \cdot P
      &= \{(0,0),(-10/21,8/21)\}
      &&\Rightarrow \{(11/21,8/21),(1,0)\}
    \\
    A^1 \cdot P
      &= \{(0,0),(-2/21,10/21)\}
      &&\Rightarrow \{(19/21,10/21),(1,0)\}
    \\
    A^2 \cdot P
      &= \{(0,0),(8/21,2/21)\}
      &&\Rightarrow \{(0,0),(8/21,2/21)\}
    \\
    A^3 \cdot P
      &= \{(0,0),(10/21,-8/21)\}
      &&\Rightarrow \{(0,1),(10/21,13/21)\}
    \\
    A^4 \cdot P
      &= \{(0,0),(2/21,-10/21)\}
      &&\Rightarrow \{(0,1),(2/21,11/21)\}
    \\
    A^5 \cdot P
      &= \{(0,0),(-8/21,-2/21)\}
      &&\Rightarrow \{(13/21,19/21),(1,1)\}
  \end{alignat*}
  Of these, $P' = \{(0,0), (8/21,2/21)\}$ is the lexicographically earliest, so it is the \emph{canonical name} for $P$, and is shown in \Cref{subfig:snubTri21}.
\end{example}

By extending our polyforms one neighboring vertex at a time and then deduplicating, we inductively construct the set of $n$-polyforms from the $(n-1)$-polyforms, beginning with the $1$-polyforms, which are precisely the vertex orbit representatives.

\begin{example}
  The three $1$-polyforms are $\{v_1\}$, $\{v_2\}$, and $\{v_3\}$, shown in \Cref{fig:snubTri1}.
  If we compute the canonical representatives of the six children of $\{v_1\}$, the three children of $\{v_2\}$, and the three children of $\{v_3\}$, we see that there are three orbits of free $2$-polyforms: $\{(0,0), (8/21,2/21)\}$, $\{(2/21,11/21),(1/3,1/3)\}$, and $\{(8/21,23/21),(13/21,19/21)\}$, which are shown in \Cref{fig:snubTri2}.

\input{figure/snub_tri2}

  Continuing in this way results in the $69$ free $5$-polyforms shown in \Cref{subfig:snubTri5}.
\end{example}

Empirically, the number of $n$-polyforms appears to be exponential in $n$, and thus the time and space complexity of the algorithm provided is exponential because it explicitly enumerates each polyform. However, in practice, we implement techniques beyond what is described above to increase the algorithmic efficiency, in part by reducing the number of objects that need to be deduplicated.
As an example of one such technique, when growing the children of a polyform, we only consider vertices that are lexicographically later than their parent.

\section*{Examples and tables}
We give two more examples of polyhedral complexes along with tables of free polyform counts for each complex.

\begin{example}
  In 1968, architect Richard Dattner designed and patented a modular play structure consisting of cuboctahedra attached along square faces.
  An example consisting of $10$ cuboctahedra is shown in \Cref{fig:DatterPlaygroundCuboctahedra}.

\input{figure/rectified_cubic}

  Richard Dattner writes, ``The simple geometry and open-endedness of PlayCubes helped to inspire new ways of playing, different interactions, and more exploration'' \cite{Dattner}.
  This is as true for the design process as it is for the play.
  In these play structures, cuboctahedra are constrained to the rectified cubic honeycomb, which contains cuboctahedral and tetrahedral cells, but Dattner's construction uses only the cuboctahedral cells, leaving octahedral voids.
  The number of $n$-polyforms containing only cuboctahedra is in direct bijection with the number of polycubes, so a $10$-celled structure consisting of only cuboctahedra could be arranged in $A038119(10) = 178\,083$ ways.
  If octahedral cells were also used, then we computed that there would be $A384254(10) = 43\,305\,326$ possible $10$-celled play structures, a more than $243$-fold increase \cite[A384254]{OEIS}. All of the $A384254(4) = 40$ possibilities for $4$-polyforms on the rectified cubic honeycomb are shown in \Cref{subfig:RectifiedCubic4forms}.
\end{example}

\begin{example}
  ZIGZAG, a modular climbing structure designed by HardBodyHang and shown in \Cref{subfig:ZIGZAG}, consists of tetrahedral and octahedral components attached together face-to-face with climbing holds on the exposed faces.
  This design was a finalist for the 2023 Design Europa Awards.
  \input{figure/octahedral-tetrahedral}

  While in practice, these climbing structures are not constrained to the tetrahedral-octahedral honeycomb, we count the $n$-celled examples that are.
  We have published the number of structures that align with the honeycomb \cite[A343909]{OEIS} and have computed up to $A343909(16) = 44\,896\,941\,979$, meaning that there are nearly $45$ billion possible ways to combine $16$ tetrahedral and octahedral cells in a way that conforms to the honeycomb, up to rotation and reflection.
\end{example}

\subsection*{Data tables}
While we have computed the initial terms of sequences for dozens of tessellations of the plane and $3$-space, in \Cref{table:data} we present terms of the sequences related to the geometries above. To see a more complete list, refer to the cross-references in OEIS sequences A385028 and A385266.

\input{figure/data_table}

\section*{Summary and Conclusions}
These polygonal and polyhedral constructions suggest other directions.
We could ask about the asymptotic growth of the number of polyforms on various tilings with respect to the number of cells, and numerous other analogs of questions that have been asked about polyominoes.
We can ask about counting polyhedral complexes that are not constrained to a periodic Euclidean graph, such as ``polytets'' which are non-intersecting regular tetrahedral complexes that have been counted by David Ellsworth \cite[A276272]{OEIS} and the analogous question about pentagonal complexes in the plane \cite{KochKurz2006}.
Other authors have considered counting objects that exist on uniform hyperbolic tilings \cite{Roldan2023}.

\subsection*{Packing Puzzles}
  Part of the increase in popularity of polyominoes was because of \textit{Pentomino packing puzzles} that appeared in Martin Gardner's columns due to Sol Golomb.
  These puzzles challenged the readers to use the $12$ free $5$-celled polyominoes to fill the $3 \times 20$, $4 \times 15$, $5 \times 12$, and $6 \times 10$ grids without gaps or overlaps.

  Since then, puzzle designers have devised related puzzles, many of which can be found under Kadon Enterprises.
  One particularly remarkable example is their soma-cube like puzzle ``Hexacube'' which uses $4$ unit cubes together with all $166$ ``one-sided hexacubes,'' which are the $6$-polyforms on the cubic honeycomb counted up to translation and rotation, but not reflection.
  This puzzle asks players to use these 170 pieces to fill in the entire $10 \times 10 \times 10$ cube with no gaps or overlaps.

  Mark Owen and Matthew Richards propose a puzzle that follows from the remarkable fact that the set of $44$ one-sided $4$-cell truncated octahedral polyforms can pack into a $3 \times 5 \times 8$ cuboid, which has precisely \(
    (3 \times 5 \times 8) + (2 \times 4 \times 7) = 44 \times 4 = 176
  \) cells.
  In their paper ``A Song of Six Splatts,'' they propose several other related packing puzzles \cite{OwenRichards1987}.
  For both truncated octahedral and rhombic dodecahedral puzzles, see Coffin's catalog \cite{Coffin2012}.

  Lastly, we discuss a puzzle with an unsolved (and perhaps unsolvable) variant idea from Andrew Kepert.
  There are $11$ \emph{one-sided} $4$-polyforms on the tetrahedral-octahedral honeycomb.
  Combined, these have a total of $20$ octahedral cells and $24$ tetrahedral cells, and they pack simply into a size-$5$ tetrahedral region, which contains $20$ octahedra and $35$ tetrahedra.
  Can they pack this region in such a way that the cells have full tetrahedral symmetry?
  This would involve filling the entire region with the exception of the tetrahedra at the four vertices, six edge-midpoints, and the center of the size-$5$ tetrahedral region. We encourage our readers to construct more puzzles that follow from the ability to enumerate all pieces of a certain size.

  Enumerating geometric structures is a natural midpoint between the polyhedral complexes that arise in art, the impulse to explore all possibilities of a system, ideas arising from geometry and group theory, and the utility of powerful modern computing.
  These suggest new possibilities in art and architecture, but also in toy and playground design.

\section*{Acknowledgments}
We thank
Pontus von Br\"omssen for his discussions and contributions to counting spherical polyforms,
along with
Jasper Bown,
Robin Houston,
Glen Whitney,
Jiangmei Wu,
the Desert Botanical Garden, and
HardBodyHang.
\bibliography{bibliography}{}
\bibliographystyle{plain}
\end{document}

%% file: figure/bisymmetric_hendecahedron.tex
\noindent
\begin{adjustbox}{margin=1em, center}
\begin{minipage}{\textwidth}
\captionsetup{type=figure}
\centering
\begin{minipage}{0.27327675\textwidth}
  \centering
  \includegraphics[width=\linewidth]{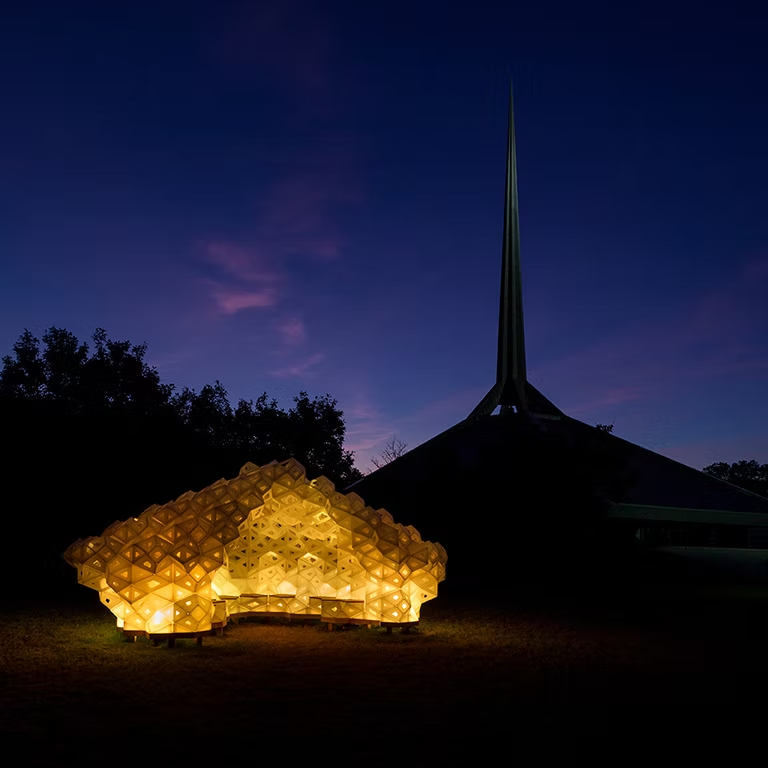}
  \subcaption{}
\end{minipage}
\qquad
\begin{minipage}{0.27327675\textwidth}
  \centering
  \includegraphics[width=\linewidth]{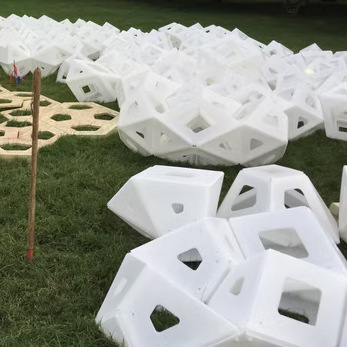}
  \subcaption{}\label{subfig:JiangmeiWu}
\end{minipage}
\qquad
\begin{minipage}{0.27327675\textwidth}
  \centering
  \includegraphics[width=\linewidth]{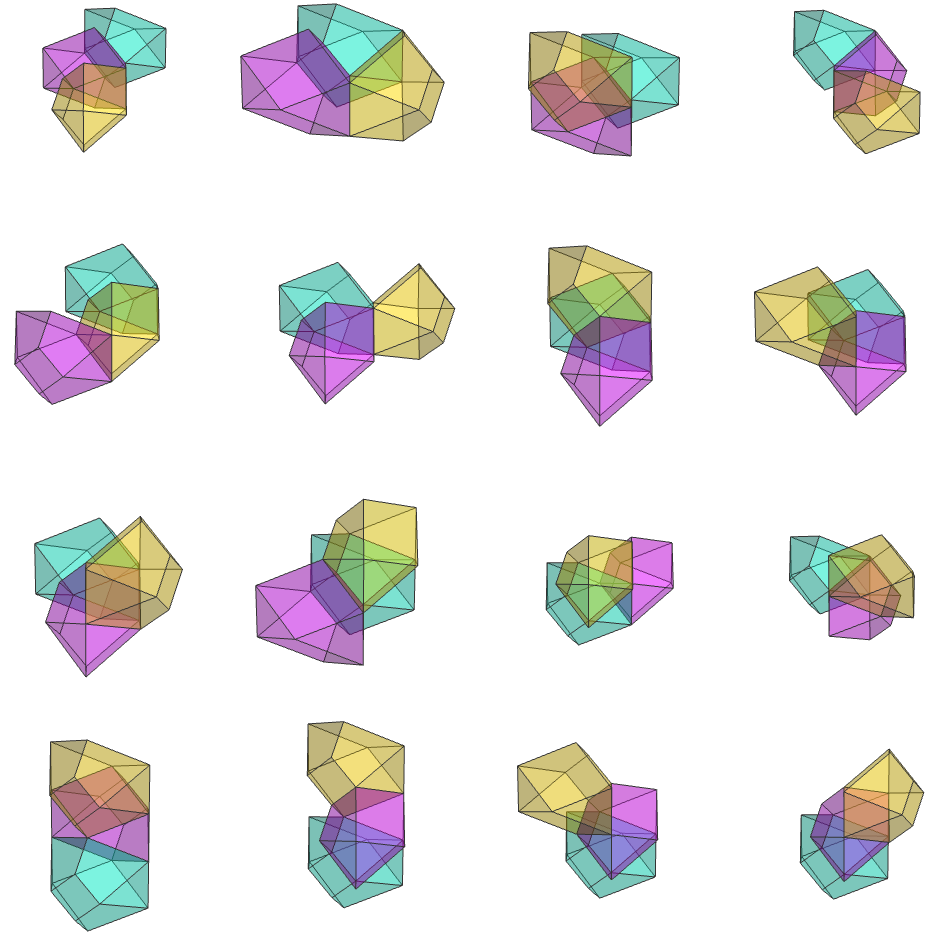}
  \subcaption{}
\end{minipage}
\caption{
  In (a), Jiangmei Wu's 2017 \emph{Synergia} consisting of hundreds of bisymmeric hendecahedra;
  in (b) a view of the hendecahedral components; and
  in (c) models of the $16$ essentially different structures we can make from three of these components. (Photographs used with permission.)
}
\label{fig:BisymmetricHendecahedron}
\end{minipage}
\end{adjustbox}

%% file: figure/truncated_octahedra.tex
\noindent
\begin{adjustbox}{margin=1em, center}
\begin{minipage}{\textwidth}
\captionsetup{type=figure}
\centering
\begin{minipage}{0.246178\textwidth}
  \centering
  \includegraphics[width=\linewidth]{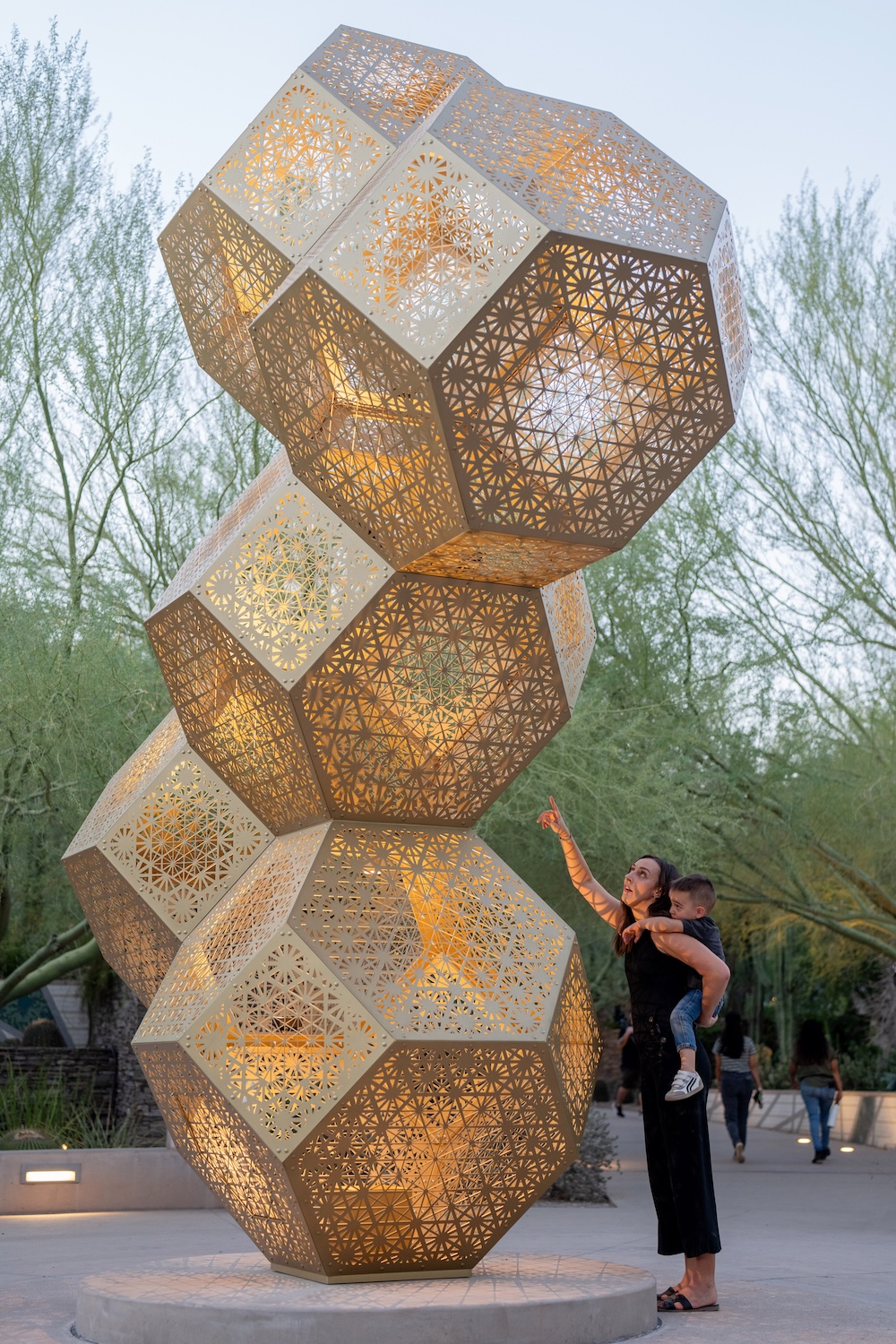}
  \subcaption{}\label{fig:TOHYBYCOZO}
\end{minipage}
\hfill
\begin{minipage}{0.369267\textwidth}
  \centering
  \includegraphics[width=\linewidth]{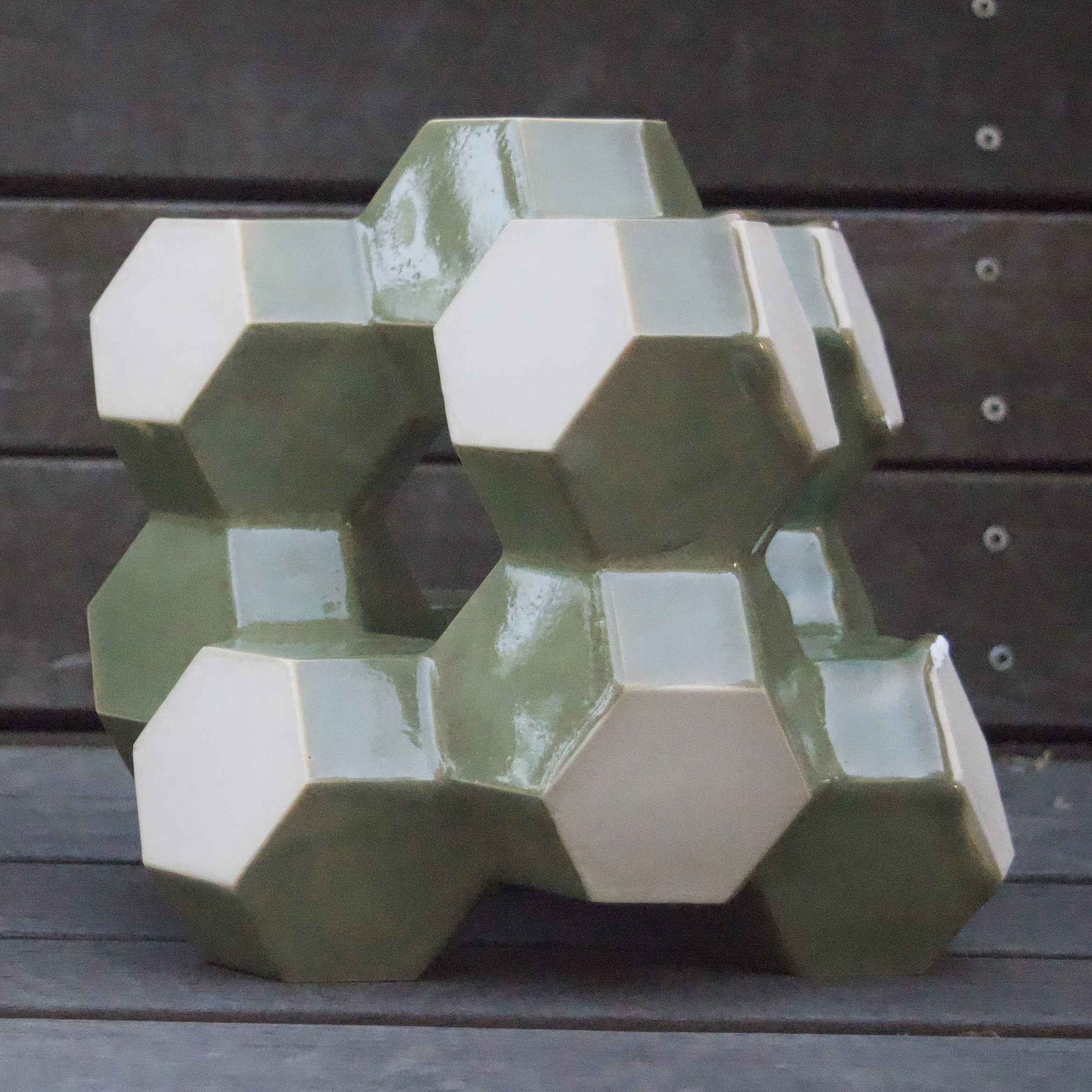}
  \subcaption{}\label{fig:TOJasper}
\end{minipage}
\hfill
\begin{minipage}{0.354555\textwidth}
  \centering
  \includegraphics[width=\linewidth]{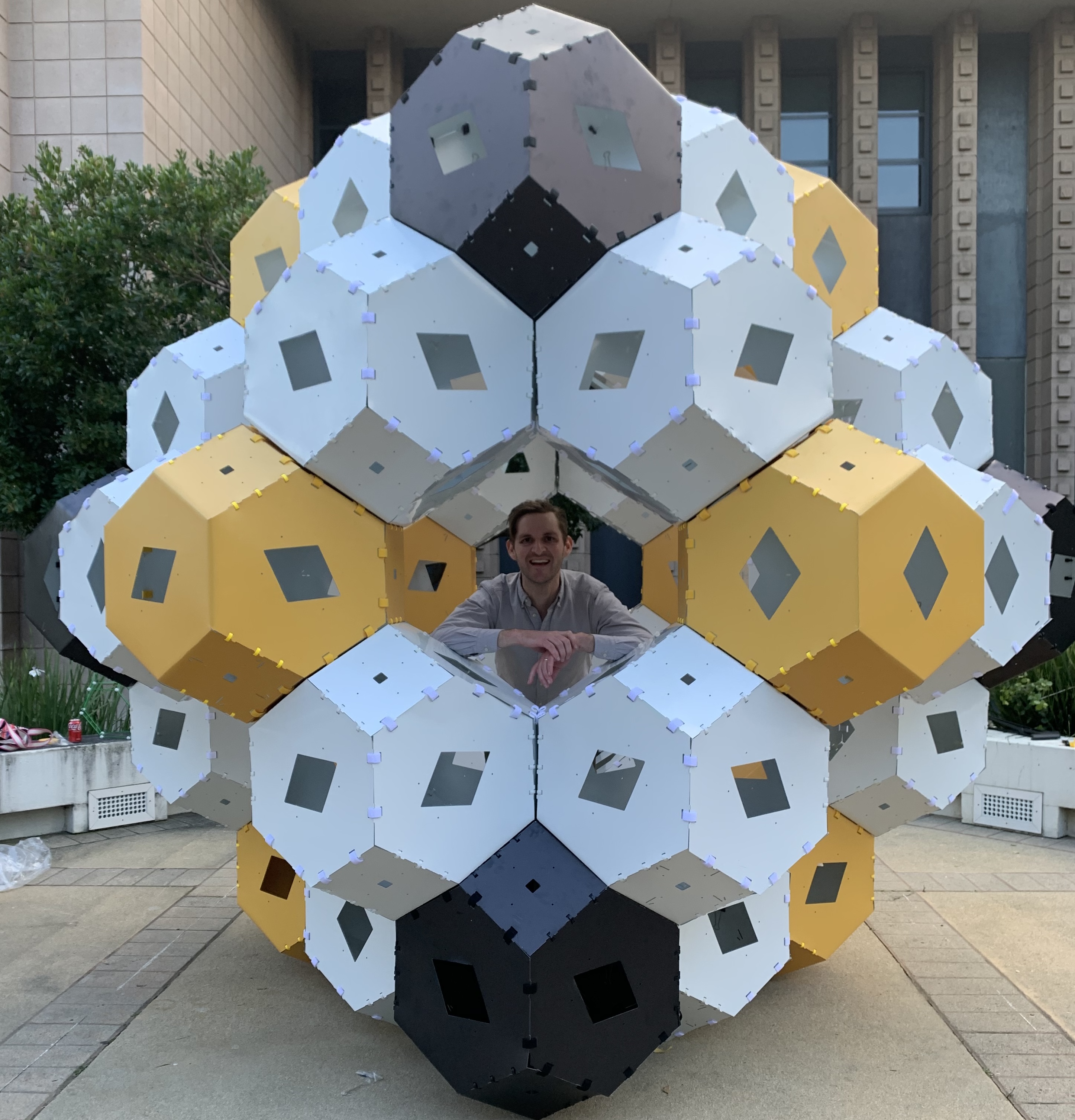}
  \subcaption{}\label{fig:TOWARD}
\end{minipage}
\caption{
  Three truncated octahedral complexes: in (a) a photograph of HYBYCOZO's sculpture of a $5$-celled truncated octahedral compound at the Desert Botanical Garden;
  in (b) a photograph of Jasper Bown's $10$-celled truncated octahedral ceramic sculpture ``Polytopes amidst a Permutahedral Tiling'' \cite{Bown2024}; and
  in
  (c) Studio Infinity's sculpture TOWARD which is a $38$-celled truncated octahedral complex.
  (Photographs used with permission.)
  \label{fig:truncatedOctahedral1}
}
\end{minipage}
\end{adjustbox}

%% file: figure/disphenoid_tetrahedron.tex
\noindent
\begin{adjustbox}{margin=1em, center}
\begin{minipage}{\textwidth}
\captionsetup{type=figure}
\centering
\begin{minipage}{0.380336\textwidth}
  \centering
  \includegraphics[width=\linewidth]{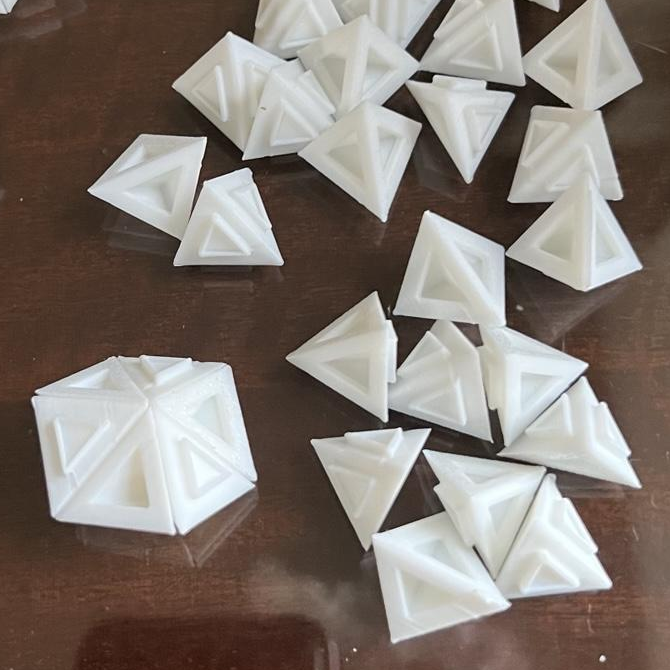}
  \subcaption{}
\end{minipage}
\hfill
\begin{minipage}{0.599664\textwidth}
  \centering
  \includegraphics[width=\linewidth]{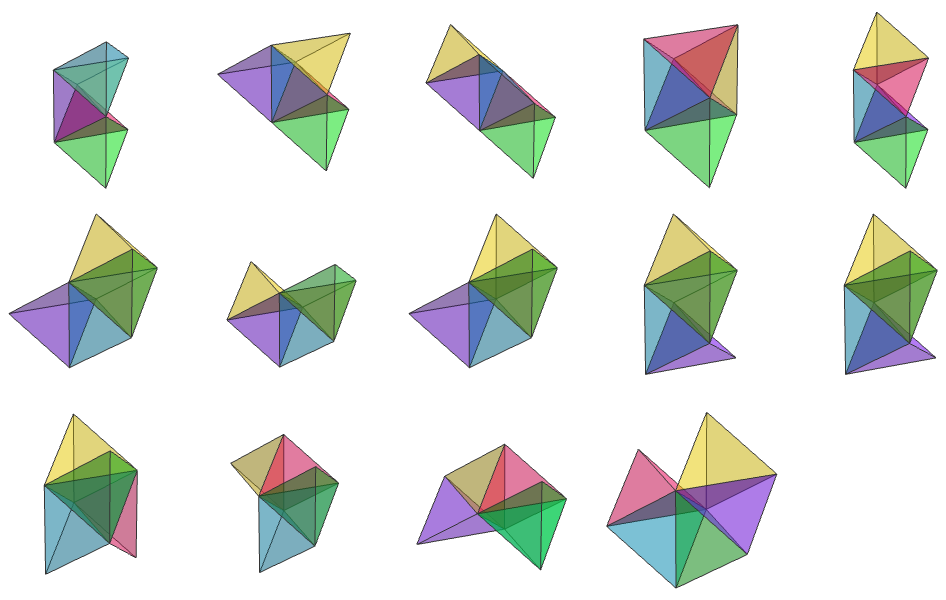}
  \subcaption{}
\end{minipage}
\caption{
  In (a) a photograph of Robin Houston's disphenoid tetrahedral building block toy \cite{Houston2025}; in (b) an illustration of the $A385024(5)=14$ constructions that can be made with $5$ blocks, up to rotation and reflection. (Photograph used with permission.)
}
\label{fig:disphenoidTetrahedral}
\end{minipage}
\end{adjustbox}

%% file: figure/snub_tri_dual_graph.tex
\noindent
\begin{adjustbox}{margin=1em, center}
\begin{minipage}{\textwidth}
\captionsetup{type=figure}
\centering
\begin{minipage}{0.409781\textwidth}
  \centering
  \includegraphics[width=\linewidth]{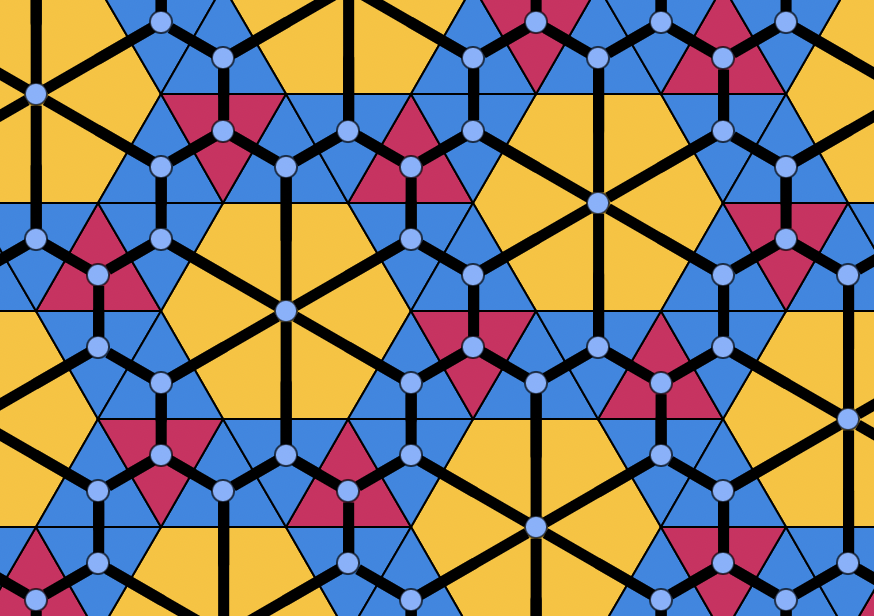}
  \subcaption{}\label{subfig:snubTrihexagonalGraph}
\end{minipage}
\hfill
\begin{minipage}{0.570219\textwidth}
  \centering
  \includegraphics[width=\linewidth]{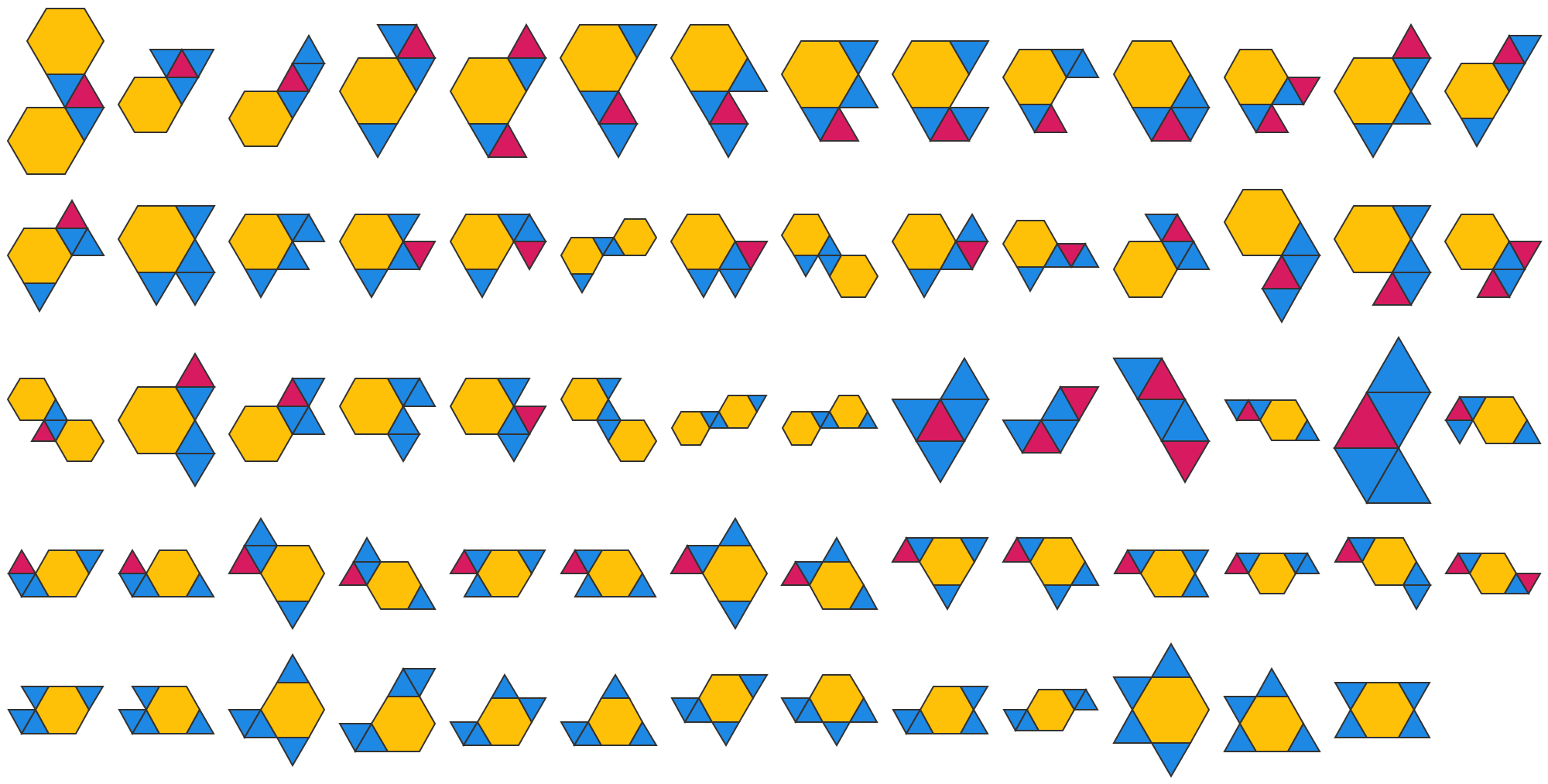}
  \subcaption{}\label{subfig:snubTri5}
\end{minipage}
\caption{
  In (a), the snub trihexagonal tiling of $\mathbb{R}^2$ overlaid with a dual graph containing a vertex at every cell and edges between adjacent cells;
  in (b) all $69$ of the $5$-polyforms on the snub trihexagonal tiling.
}
\label{fig:SnubTrihexagonal}
\end{minipage}
\end{adjustbox}

%% file: figure/snub_tri1.tex
\noindent
\begin{adjustbox}{margin=1em, center}
\begin{minipage}{\textwidth}
\captionsetup{type=figure}
\centering
~
\hfill
\begin{minipage}{0.20\textwidth}
  \centering
  \includegraphics[width=\linewidth]{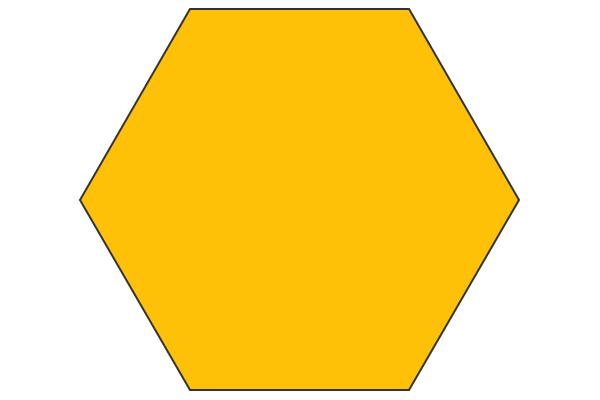}
  \subcaption{
    \(v_1 = \left(0,0\right)\)
  }\label{fig:snubTri11}
\end{minipage}
\hfill
\begin{minipage}{0.20\textwidth}
  \centering
  \includegraphics[width=\linewidth]{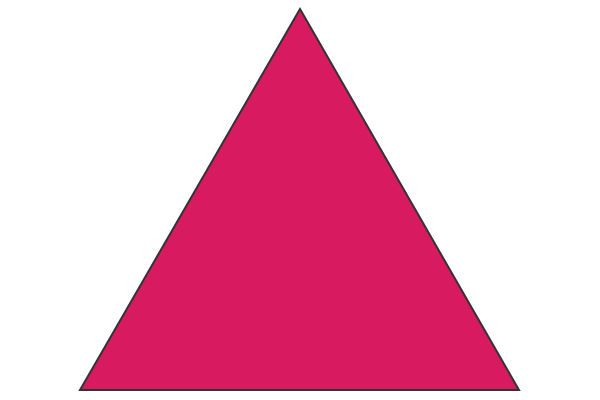}
  \subcaption{
    \(v_2 = \left(\frac{1}{3},\frac{1}{3}\right)\)
  }\label{fig:snubTri12}
\end{minipage}
\hfill
\begin{minipage}{0.20\textwidth}
  \centering
  \includegraphics[width=\linewidth]{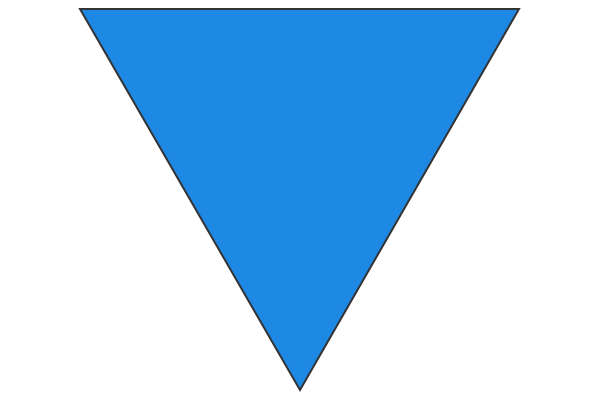}
  \subcaption{
    \(v_3 = \left(\frac{2}{21},\frac{11}{21}\right)\)
  }\label{fig:snubTri13}
\end{minipage}
\hfill
~
\caption{
  The three distinct cells that correspond to the three orbits of vertices in the dual graph under its automorphism group. These are (a) hexagonal cells, (b) triangular cells with no hexagonal neighbors, and (c) triangular cells with hexagonal neighbors.
}
\label{fig:snubTri1}
\end{minipage}
\end{adjustbox}

%% file: figure/snub_tri2.tex
\noindent
\begin{adjustbox}{margin=1em, center}
\begin{minipage}{\textwidth}
\captionsetup{type=figure}
\centering
~
\hfill
\begin{minipage}{0.30\textwidth}
  \centering
  \includegraphics[width=0.66\linewidth]{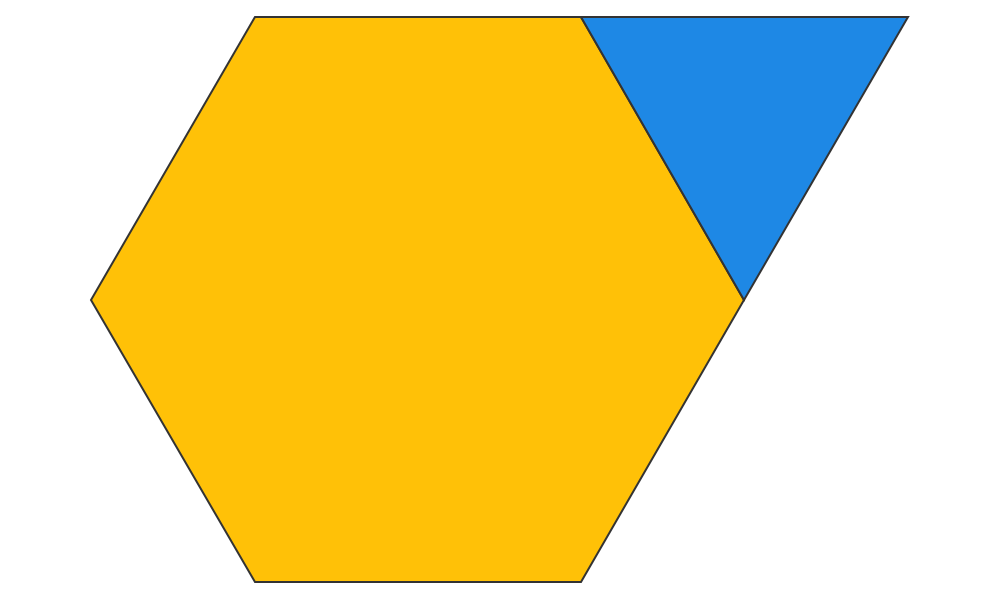}
  \subcaption{$\{(0,0), (8/21,2/21)\}$}
  \label{subfig:snubTri21}
\end{minipage}
\hfill
\begin{minipage}{0.30\textwidth}
  \centering
  \includegraphics[width=0.66\linewidth]{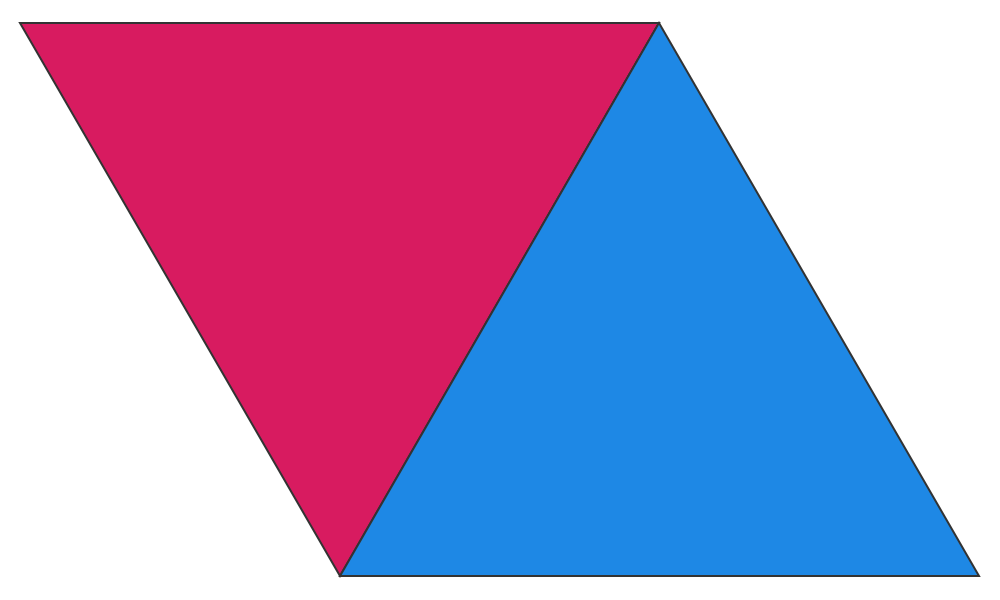}
  \subcaption{$\{(2/21,11/21),(1/3,1/3)\}$}
  \label{subfig:snubTri22}
\end{minipage}
\hfill
\begin{minipage}{0.30\textwidth}
  \centering
  \includegraphics[width=0.66\linewidth]{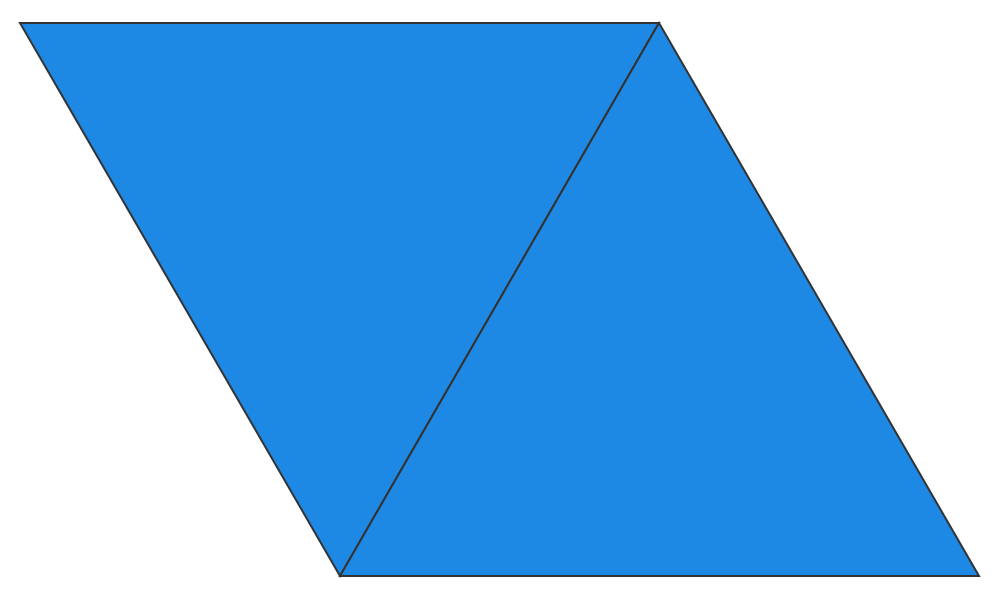}
  \subcaption{$\{(8/21,23/21),(13/21,19/21)\}$}
  \label{subfig:snubTri23}
\end{minipage}
\hfill
~
\caption{
  The three free $2$-polyforms on the snub trihexagonal tiling with their canonical names.
}
\label{fig:snubTri2}
\end{minipage}
\end{adjustbox}

%% file: figure/rectified_cubic.tex
\noindent
\begin{adjustbox}{margin=1em, center}
\begin{minipage}{\textwidth}
  \captionsetup{type=figure}
  \centering
  \begin{minipage}{0.277143\textwidth}
    \includegraphics[width=\textwidth]{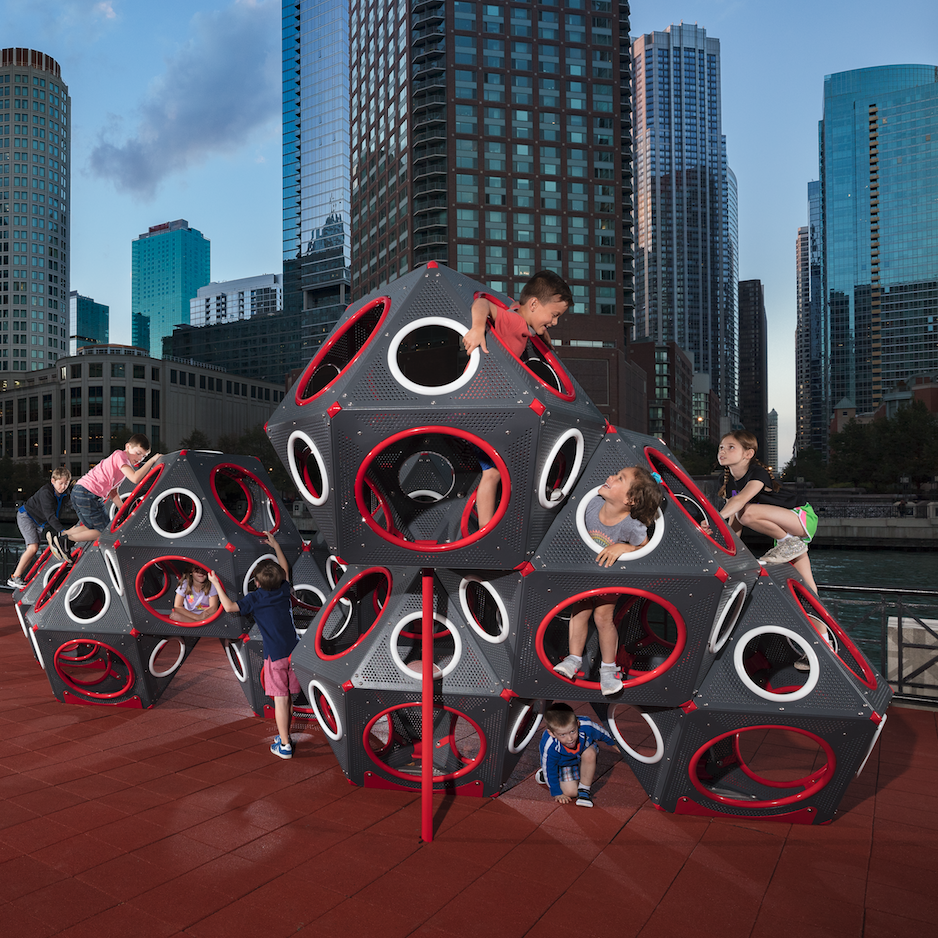}
    \subcaption{}
    \label{subfig:DattnerPlayground}
  \end{minipage}
  \begin{minipage}{0.692857\textwidth}
    \includegraphics[width=\textwidth]{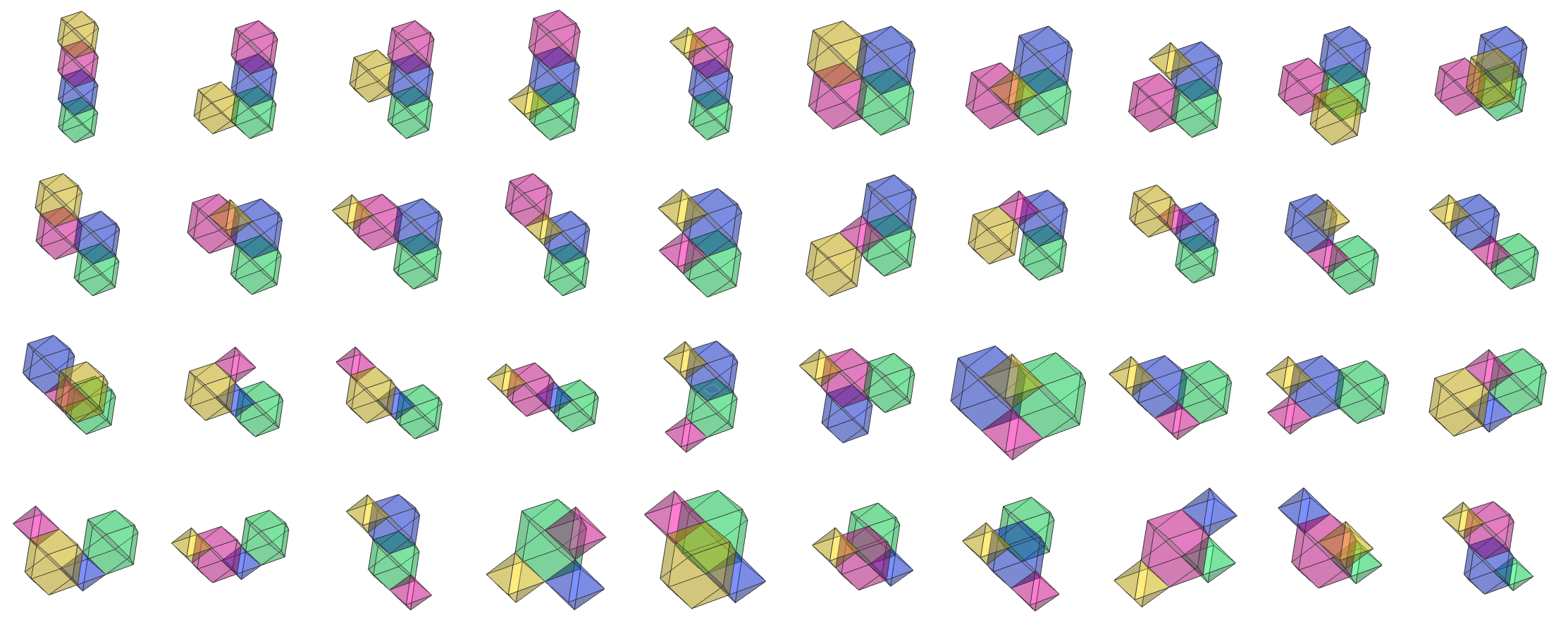}
    \subcaption{}
    \label{subfig:RectifiedCubic4forms}
  \end{minipage}
  \caption{
    A play structure in Chicago called PlayCubes (a) consisting of $10$ cuboctahedra connected face-to-face based on a 1968 design by architect Richard Dattner and manufactured by Playworld; and an illustration (b) of the $A384254(4)=40$ free $4$-polyforms on the rectified cubic (cuboctahedral-octahedral) honeycomb.
  }
  \label{fig:DatterPlaygroundCuboctahedra}
\end{minipage}
\end{adjustbox}

%% file: figure/octahedral-tetrahedral.tex
\noindent
\begin{adjustbox}{margin=1em, center}
\begin{minipage}{\textwidth}
  \captionsetup{type=figure}
  \centering
  \begin{minipage}[t]{0.364051\textwidth}
      \includegraphics[width=\textwidth]{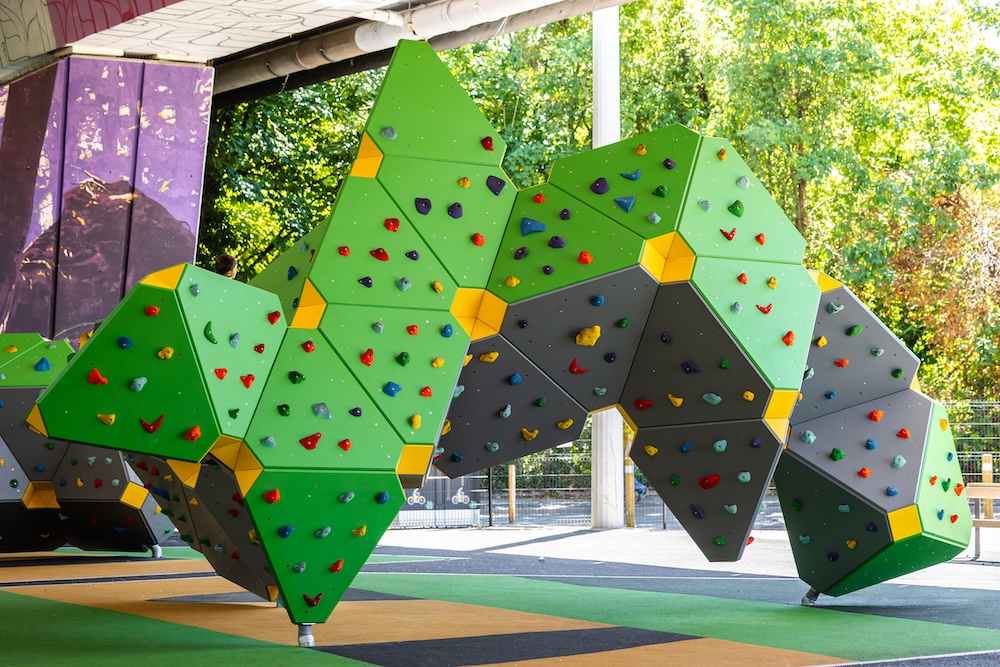}
      \subcaption{}
      \label{subfig:ZIGZAG}
  \end{minipage}
  \begin{minipage}[t]{0.605949\textwidth}
      \includegraphics[width=\textwidth]{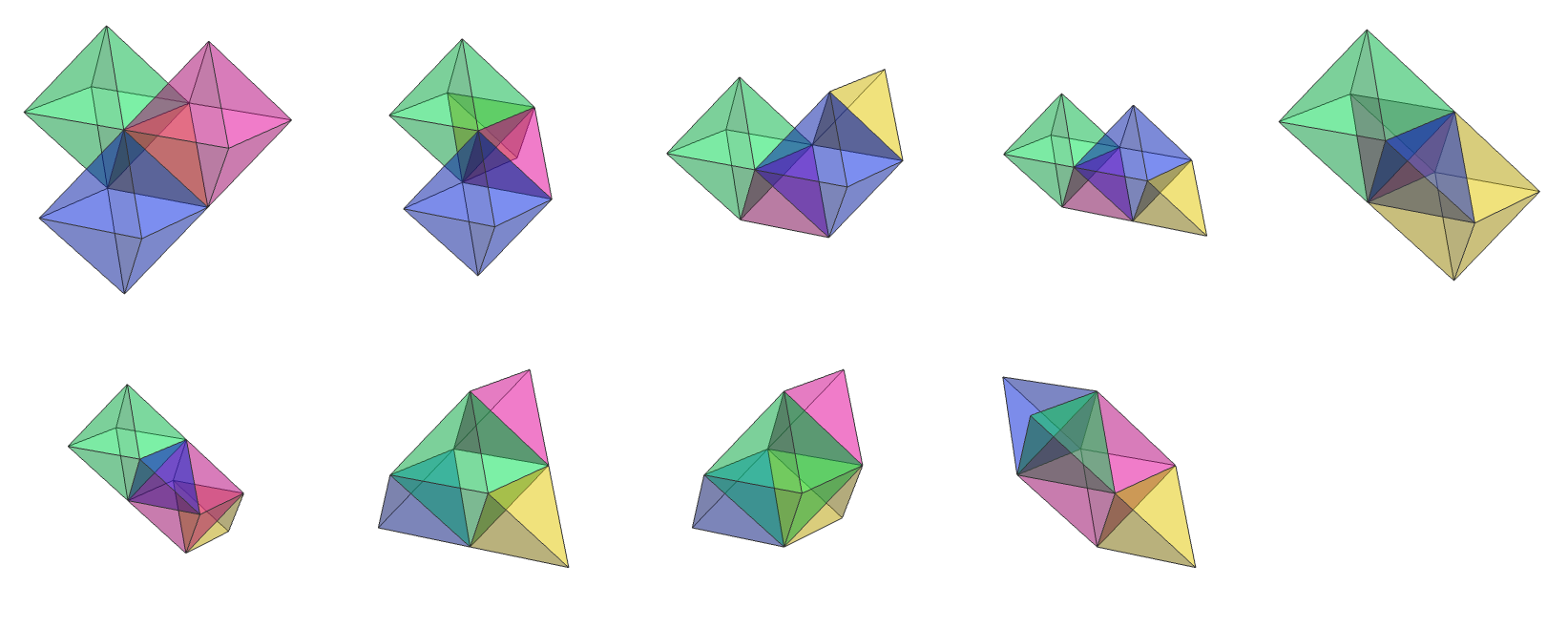}
      \subcaption{}
      \label{subfig:octtet4}
  \end{minipage}
  \caption{
    In (a) ZIGZAG by HardBodyHang is a climbing structure composed of octahedral and tetrahedral modules. (In practice, these structures are not octahedral-tetrahedral honeycomb.) In (b) the $A343909(4) = 9$ free $4$-polyforms on the tetrahedral-octahedral honeycomb. (Photograph used with permission.)
  }
  \label{fig:octtet}
\end{minipage}
\end{adjustbox}

%% file: figure/data_table.tex
\noindent
\begin{adjustbox}{margin=1em, center}
\begin{minipage}{\textwidth}
  \captionsetup{type=table}
  \caption{The initial terms of several sequences giving the number of free $n$-polyforms (for $n \geq 1$).}
  \label{table:data}
  \centering
  \begin{tabular}{lll|rrrrrrr}
      \hline
     & & &  \multicolumn{7}{c}{$n$-polyforms} \\
    Tessellation & OEIS & Figure & 1 & 2 & 3 & 4 & 5 & 6 & 7 \\
    \hline
    Bisymmetric hendecahedral
      & \href{https://oeis.org/A385028}{A385028}
      & \Cref{fig:BisymmetricHendecahedron}
      & 1 & 4 & 16 & 116 & 903 & 8551 & 85365
      \\
    Truncated octahedral
      & \href{https://oeis.org/A038181}{A038181}
      & \Cref{fig:truncatedOctahedral1}
      & 1 & 2 & 6 & 35 & 251 & 2602 & 30900
      % 30900, 400818
      \\
    Disphenoid tetrahedral
      & \href{https://oeis.org/A385024}{A385024}
      & \Cref{fig:disphenoidTetrahedral}
      & 1 & 1 & 2 & 5 & 14 & 47 & 172
      % , 691, 2881, 12449, 54782, 244992
      \\
    Snub trihexagonal
      & \href{https://oeis.org/A383908}{A383908}
      & \Cref{fig:SnubTrihexagonal}
      & 3 & 3 & 7 & 23 & 69 & 228 & 766
      % , 2642, 9309, 33382, 120629
      \\
    Rectified cubic
      & \href{https://oeis.org/A384254}{A384254}
      & \Cref{fig:DatterPlaygroundCuboctahedra}
      & 2 & 2 & 9 & 40 & 290 & 2529 & 26629
      % 301289, 3568048, 43305326, 534671742, 6684869463
      \\
    Cubic
      & \href{https://oeis.org/A038119}{A038119}
      & \Cref{fig:DatterPlaygroundCuboctahedra}
      & 1 & 1 & 2 & 7 & 23 & 112 & 607
      % 3811, 25413, 178083
      \\
    Tetrahedral-octahedral
      & \href{https://oeis.org/A343909}{A343909}
      & \Cref{fig:octtet}
      & 2 & 1 & 4 & 9 & 44 & 195 & 1186
      % , 7385, 49444, 337504, 2353664
      \\
    \hline
  \end{tabular}
\end{minipage}
\end{adjustbox}